\documentclass[11pt, a4paper]{amsart} 
\usepackage[utf8]{inputenc}
\usepackage[T1]{fontenc}
\usepackage{listings}
\usepackage{xcolor}
\usepackage{caption}
\usepackage{tikz}
\usetikzlibrary{er,positioning}
\usetikzlibrary{calc,intersections,through,backgrounds}

\definecolor{codegreen}{rgb}{0,0.6,0}
\definecolor{codegray}{rgb}{0.5,0.5,0.5}
\definecolor{codepurple}{rgb}{0.58,0,0.82}
\definecolor{backcolour}{rgb}{0.95,0.95,0.92}
 
\lstdefinestyle{mystyle}{
    backgroundcolor=\color{backcolour},   
    commentstyle=\color{codegreen},
    keywordstyle=\color{magenta},
    numberstyle=\tiny\color{codegray},
    stringstyle=\color{codepurple},
    basicstyle= \tiny,  
    breakatwhitespace=false,         
    breaklines=true,                 
    captionpos=b,                    
    keepspaces=true,                 
    numbers=left,                    
    numbersep=5pt,                  
    showspaces=false,                
    showstringspaces=false,
    showtabs=false,                  
    tabsize=1
}
 
\usepackage[utf8]{inputenc}
\usepackage[T1]{fontenc}
\usepackage{listings}
\usepackage{amssymb,epsfig}
\usepackage{verbatim}
\usepackage{fancyvrb}

\usepackage[T1]{fontenc}
\usepackage{epsfig,multicol, graphics}
\usepackage{hyperref}
\usepackage{pdfpages}
\usepackage{marginnote}
\input xypic
\xyoption{all}
\usepackage{graphics}
\usepackage{graphicx}
\usepackage{color}

\newtheorem{thm}{Theorem}
\newtheorem{pro}[thm]{Proposition}

\newtheorem{ex}[thm]{Example}
\newtheorem{exe}[thm]{Example}

\newcommand{\F}{{\mathcal F}}

\newcommand{\x}{{\bf x}}

\newcommand{\plano}{{\mathbb P}^2}

\renewcommand{\P}{{\mathbb P}}
\renewcommand{\P}{{\mathbb P}}
\newcommand{\tor}{\xymatrix{\ar@{-->}[r]&}}
 
\definecolor{grass}{rgb}{0.14,0.72,0.2}

\begin{document}

\lstset{language=Python}

{\tiny \title{Effective Integrability of Lins Neto's Family  of   Foliations}}
\author{by\,   Lu\'is Gustavo Mendes and Liliana Puchuri}

\begin{abstract}
	A. Lins Neto presented in \cite{LinsNeto}   a  $1$-dimensional family  of degree four foliations on the complex projective plane $\F_{t \in \overline{\mathbb{C}}}$ with non-degenerate singularities of fixed analytic type,   whose set of parameters $t$ for which $\F_t$ is an elliptic pencil  is dense and countable. In \cite{Mc} and \cite{Guillot}, M. McQuillan and A. Guillot showed that the family lifts  to  linear foliations  on the abelian surface $E \times E$, where $E = \mathbb{C}/\Gamma$, $\Gamma = < 1 , \tau>$  and  $\tau$ is  a primitive 3rd root of unity,  the parameters  for which $\F_t$ are elliptic pencils   being  $t\in \mathbb{Q}(\tau) \cup {\infty}$. In  \cite{Puchuri}, the second author   gave  a closed  formula for the degree  of the elliptic curves of $\F_t$  a function of $t \in  \mathbb{Q}(\tau)$.  In this work we   determine    degree,    positions and  multiplicities of singularities of the elliptic curves of $\F_ t$,   for any given $t \in \mathbb{Z}(\tau)$   in algorithmical way  implemented in  \emph{Python}. And also we obtain  the  explicit expressions for the generators of the elliptic pencils, using the  \emph{Singular} software. Our  constructions depend on the effect of quadratic Cremona maps on  the family of foliations  $\F_t$. 
	
\end{abstract}

\maketitle

\section{Introduction and results}

An important consequence of the $1$-parameter  family  of degree four foliations $\F_t$  on the projective plane found in \cite{LinsNeto}  by A. Lins Neto was to show that the local analytic types of  singularities do not give enough   information on the existence of global first integrals neither  on  the degrees of possible  invariant algebraic  curves. This was a negative answer to a problem studied  by  H. Poincaré in two papers published in \emph{Rendiconti del Circolo Matematico di Palermo} (\cite{Poincare1}, \cite{Poincare2}).

Just for a dense and countable set of parameters  $D \subset \overline{\mathbb{C}}$ the foliations $\F_{t\in  D}$  have  global   first integrals. In fact, the foliations $\F_{t\in D}$  are elliptic pencils, whose generic elements have unbounded degrees as $t$ varies in $D$.

As remarked in \cite{Mc} and \cite{Guillot}  by M. McQuillan and A. Guillot,  the Lins Neto's family of  degree  four foliations $\F_t$  of $\plano$  lifts to a family of linear foliations on the abelian surface  $E \times E$, where $E = \mathbb{C}/\Gamma$, $\Gamma = < 1 , \tau>$  and  $\tau$ is the  primitive 3rd root of unity. And the parameters $t$ corresponding  to  elliptic pencils  were  described as $t \in \mathbb{Q}(\tau) \cup \infty$.

In \cite{Puchuri}, there is a closed formula for the degree   of the generic element of the elliptic pencil  $\F_ t$  in terms of the arithmetic of $t \in \mathbb{Q}(\tau)$.

The question we consider here is, for $t  = m+ n \tau  \in \mathbb{Z}(\tau)$, to determine in  algorithmical  way the  degree, positions and multiplicities of singularities of the generic elements of the elliptic  pencil $\F_{m+ n \tau}$,  and even to give the  expression of the generators  of such  elliptic pencil.

Our main  result is:

\begin{thm} 

Let   $t   = m + n \tau \in   \mathbb{Z} (\tau)$, where  $\tau:= e^{\frac{2 \pi I}{3}}$.  Starting  with the  elliptic  pencils  of cubics 
{\small \[\begin{cases} \F_1 :   c_1 \cdot (y-x)(z-x) (y-z) + c_2 \cdot (y-\tau x)(z-\tau^2 x) (z-\tau y) = 0 \\
\F_\infty: c_1 \cdot (y-x) (y-\tau x) (y-\tau^2 x) + c_2 \cdot (y-z)(y-\tau z) (y- \tau^2 z)=0, \quad (c_1:c_2) \in \overline{\mathbb{C}},\end{cases} \]}
and transforming  them by  a  number of sucessive applications of  the   quadratic Cremona maps  $Q_i$, $i \in     \{ 1, \tau, \tau^2, \infty\}$,  
\[\begin{cases}   Q_1(x:y:z) = (y^2 -x z : x^2 - y z: z^2- x y  ) \\ Q_\tau (x:y:z)= (\tau y^2 - x z, \tau x^2 - y z: z^2- \tau^2 x y    ) \\ Q_{\tau^2} (x:y:z)= (\tau^2 y^2 - x z, \tau^2 x^2 - yz : z^2- \tau x y) \\ Q_\infty(x:y:z)=  ( y z: x z: x y )   \end{cases}   \]
in an   algorithmicaly determined order,   we obtain the degree, position and multiplicities of the generic element  of  the  elliptic pencil $\F_{m+n \tau}$, as well the  expression of   its  generator curves, that is, the  first integral of the foliation.

\end{thm}\label{principal}

For example, the generic element of the elliptic pencil $\F_{-40 + 160 \tau}$ has degree  $100806$ and its singular set is given by  three ordinary $33841$-uple points, three ordinary $33241$-uple points,  three ordinary $33721$-points   and three ordinary  $3$-uple points. The pencil $\F_{-40 + 160 \tau}$ is the strict transform  after $241$ applications of the quadratic maps  $Q_i$ maps to the pencil of cubics $\F_\infty$.

Another example,  the generic element of the elliptic pencil $\F_{180 - 110 \tau}$ has degree  $ 64302 $ and its singular set is given by  three ordinary $21277$-uple points, three ordinary $21567$-uple points and   three ordinary $21457$-points. The pencil $\F_{180 - 110 \tau}$ is the strict transform  after $312$ applications of the quadratics $Q_i$ maps to the pencil of cubics $\F_1$. 

Section \ref{sectionexamples} gives more   examples and present the algorithm of Theorem \ref{principal}  implemented in \emph{Python}, with which the reader can carry out an indefinite number of examples and experiments. For getting the expression of generator curves, i.e. first integrals, the implementation is done in   \emph{Singular} software.

\emph{Question:} At the end of Section \ref{taut} we explain  that our method of realization of the elliptic pencils $\F_t$ extends  to an infinite  number of parameters of the form $t= p + q \tau \in \mathbb{Q}(\tau)$. But the  question that  remains is  how to reach in algorithmical way the elliptic pencils corresponding to  all parameters $t \in \mathbb{Q}(\tau)$. 

In a second work \cite{MPNegative} we give applications of our method to the study of negative curves on rational surfaces.

{\emph{Acknowledgments:}  The first  author  thanks the  participants of the online Seminar  \emph{Painlev\'e - Stockholm - CNRS}  for an  invitation  to speak  on  Painlev\'e 's results on pencils of curves, which gave rise  to the  the question  studied here. Special thanks to Adolfo Guillot, for a letter on the lift of the foliations to abelian surfaces and on degrees of the first integrals. He also thanks  Vitalino Cesca Filho for help on the  \emph{Singular} software. We both  thank Orestes Bueno for improving the quality of figures. }

{\small \tableofcontents }

\section{Background material on foliations}\label{background}

\quad

A singular holomorphic foliation $\F$ of  a smooth   projective surface $M$  can  given by a finite  open covering $\{U_i\}$ and local differential equations \[\omega_i(x_i,y_i) = a_i(x_i,y_i) dx_i + b_i(x_i,y_i) dy_i  = 0, \quad  a_i,b_i  \in  \mathcal{O}(U_i),\, \mbox{with}\,  gcd(a_i,b_i)=1,\]  such that along  $U_i \cap U_j \neq \emptyset $  it holds $\omega_i = g_{ij} \, \omega_i$   for  $g_{ij}\in \mathcal{O}^*(U_i\cap U_j)$.  The conditions   $gcd(a_i,b_i)=1$ assure that the  singular set  of the foliation  $\mbox{Sing}(\F) \cap U_i = \{ \omega_i(p) = 0\}$   is finite. 

We  pass from a local  $1$-form $\omega(x,y) = a(x,y) dx + b(x,y) dy $  to its dual vector field $v = b(x,y) \,  \frac{\partial }{\partial x}  - a(x,y) \frac{\partial }{\partial y}$.  And define a  singularity $p$   of  the  foliation  as being of    \emph{reduced type}  if  $v$  has  non-trivial linear part,   with  at least one non-zero eigenvalue $\lambda_1$ and $\frac{\lambda_2}{\lambda_1}\not\in \mathbb{Q}^+$; and  of  \emph{reduced nondegenerate type} if it is reduced and both   eigenvalues are not zero. For instance, in this paper we shall deal with reduced nondegenerate  singularities of type $\omega = 3 y dx + x dy + h.o.t = 0$, which have local holomorphic first integrals of type $x^3 y = c$. 

After a finite number of blowing ups,  any singularity of the foliation is replaced by  a number of reduced singularities along the  exceptional divisor (Seidenberg's reduction of singularities). In this paper, only reduced nondegenerate points appear after a reduction of singularities. 

By \emph{dicritical} we mean a  singularity of foliation whose blow up $\sigma$  produces a non-invariant exceptional line $E= \sigma^{-1}(p)$. For instance, a     \emph{radial point} $\omega = y dx - x dy + h.o.t. =0$ is a non-reduced singularity whose blow up   gives rise to a  foliation (with isolated singularities)  completely transversal to the exceptional line. Radial points have local meromorphic first integrals of type $\frac{x}{y}= c$.

Let $\nu(p)\geq 0$ be  the order of the first non-trival jet of  a $1$-form $\omega = a(x,y)  dx + b(x,y)  dx$ defining the foliation $\F$ around $p$ by $  \omega =0$.   
Define $l(p,\F) := \nu(p)$ if $p$ is  not dicritical or   $l(p,\F) := \nu(p) +1$ if $p$ is dicritical. For example, reduced singular points have $\nu(p,\F)= l(p,\F) =1$;  radial point  has  $\nu(p,\F)=1$ and  $l(\F,p)=2$. In this work we also deal   with  non-reduced points of local  type $\omega = d(x  y \, (x-y) ) = 0$,  for  which  $\nu(p,\F)= l(p,\F) =2$.

The \emph{multiplicity} (or \emph{Milnor number})  $\mu(p,\F)$ of a singularity $p$  of a  foliation  $\F:\, \omega = a(x,y) dx + b(x,y) dy =0$  is the multiplicity of intersection of the curves $a(x,y) =0$ and $b(x,y)=0$ at $p$.  

The \emph{degree} $\mbox{deg}(\F)$ of a foliation on the complex projective plane is the number of tangencies of $\F$ and a generic projective line.

We recall:

\begin{pro}\label{Darboux} (Darboux's formula for foliations of the plane) For $\F$ be a   singular holomorphic foliation  of $\plano$ (with finite singular set),  
{\small \[ \mbox{deg}^2(\F) + \mbox{deg}(\F) + 1 = \sum_{p \in sing(\F)}
   \mu(p,\F) \]}
\end{pro}

The general foliations we shall encounter in this work have degree four,   twelve radial points and nine points of type  $\omega = 3 y dx + x dy + h.o.t=0$.

Next result shall be examplified   in Section  \ref{ellipticpencilssurfaces}:  it determines  the degree \emph{as a foliation} of a pencil of plane curves and is associated to G. Darboux (Bull. Sc. Math. 1876). It shows  that the presence of multiple components in curves of the pencil   lowers the degree of the pencil \emph{as a foliation}:

\begin{pro}\label{darboux}  (Darboux' formula for pencils)  Let  
\[c_1 \cdot F(x:y:z) + c_2 \cdot   G(x:y:z) = 0,  \quad (c_1:c_2) \in \overline{\mathbb{C}}\]
be a pencil of curves, whose   generic element is an  irreducible curve  of degree $\mbox{deg}(C)= \mbox{deg}(F) = \mbox{deg}(G)$. 
Let $C_s =\sum_k \alpha_ {s,k} \cdot  C_{s,k} $ be a decomposition in irreducible factors of special elements $C_s$, i.e. $C_s$ is an element of the pencil  having  some isolated  singularity or some  multiple component (i.e. $\alpha_{s,k} \geq 2$). 
Then the degree as foliation of this pencil of curves is
\[\mbox{deg}(\mathcal{F}) = 2 \cdot \mbox{deg}(C) - 2 -   \sum_{s,k}  ( \alpha_{s,k} -1) \cdot \mbox{deg}(C_{s,k}) \] 
where the sum runs over all special elements. 
\end{pro}

For instance: a  generic pencil of cubics has, as a foliation, $\mbox{deg}(\F)= 2\cdot 3 - 2 =4$. But the pencil $c_1 \cdot x  y  z + c_2 \cdot z^3 = 0$ has, as a foliation, degree $2\cdot 3 -2 - 2 =1$, thanks to the $3$-uple line $z=0$.

\section{The dual Hesse arrangement}\label{arrangement}

The  \emph{dual Hesse arrangement} of projective lines and points   on the complex projective plane 
is composed by a set of nine    lines $\mathcal{L}_9$  intersecting at  twelve  points $\mathcal{P}_{12}$. Each one of the twelve points  is a triple point of  the set of  lines  and  on  each line there are four of the twelve points. In the research field of arrangements these incidences are   usualy denoted  by $(12_3, 9_4  )$.

The   lines $\mathcal{L}_9$  can be  given  in homogeneous coordinates $(x:y:z)$ by 
\[(\star)\begin{cases}
l_1:= y-x =0;\quad  l_2:= y-\tau \cdot  x =0;\quad   l_3:= y-\tau^2 \cdot x =0 \\ 
  m_1:= z-x =0;\,\,  m_2:= z- \tau \cdot x =0; \,\, m_3:= z-\tau^2 \cdot x=0\\
 n_1:= z-y=0;\quad  n_2:= z-\tau \cdot y =0;\,\,  n_3:= z-\tau^2 \cdot y =0;
\end{cases}\]~
where $\tau = e^{\frac{2 \pi I}{3}}$, and an  equation for all  the set $\mathcal{L}_9$ is 
\[(x^3 - z^3) \cdot (y^3-z^3) \cdot (x^3-y^3) = 0 \]
 
The twelve points $\mathcal{P}_{12}$ can be decomposed  in four sets of three points 
\[\mathcal{P}_{12} = \mathcal{P}_3 (1)\cup \mathcal{P}_3 (\tau) \cup \mathcal{P}_{3}(\tau^2) \cup \mathcal{P}_{3} (\infty)\]  
as follows: 
\[(\star\star)
\begin{cases} 
\mathcal{P}_3(1):\,  (1:1:1)= l_1 \cap m_1\cap n_1,\quad  (1:\tau:\tau^2)= l_2\cap m_3\cap n_2,\quad  (1:\tau^2:\tau) = l_3\cap m_2\cap n_3 \\ 
\mathcal{P}_3 (\tau):\,  (1:\tau:1)= l_2\cap m_1\cap n_3 ,\quad  (1:1:\tau)= l_1\cap m_2\cap n_2 ,\quad    (\tau:1:1)=  l_3 \cap m_3 \cap n_1   \\
\mathcal{P}_3 (\tau^2):\, (1:1:\tau^2)= l_1\cap m_3 \cap n_3,\,\,  (\tau^2:1: 1) =  l_2 \cap m_2 \cap n_1 \,\,\, \,  \,   (1:\tau^2:1) = l_3 \cap m_1 \cap n_2 \\
\mathcal{P}_3 (\infty):\,  (0:0:1) = \, l_1 \cap l_2 \cap l_3,\quad \,   (0:1:0)= m_1 \cap m_2 \cap m_3   ,\,\,    (1:0:0)= n_1 \cap n_2 \cap n_3  
\end{cases}
\]

The dual  Hesse arrangement is projectively rigid (cf.  \cite{LinsNeto} Prop. 1 or  \cite{Magdalena} Th. 1): that is,  the incidences  $(12_3, 9_4  )$ determine the dual Hesse arrangement up to automorphism of the complex projective plane.

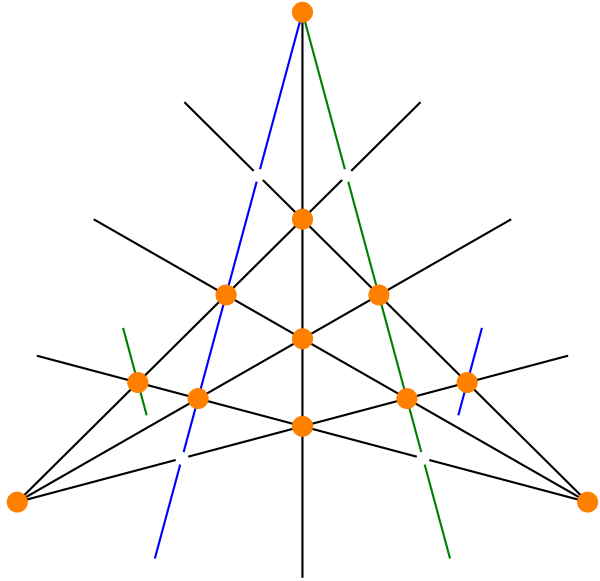
\begin{figure}
\begin{tikzpicture}[scale=1.5]
	\coordinate (O) at (0,0);
	\draw[thick, name path=O15] (O)-- (15:5);	
	\draw[thick, name path=O30] (O) -- (30:5);
	\draw[thick, name path=O45] (O) -- (45:5);
	\coordinate (V1) at (5,0);	
	\draw[thick, name path=V115] (V1) -- +(135:5);
	\draw[thick, name path=V130] (V1) -- +(150:5);
	\draw[thick, name path=V145] (V1) -- +(165:5);
	\coordinate (V2) at (60:5cm);
	\draw[thick,blue, name path=V215] (V2) -- +(255:5);
	\draw[thick, name path=V230] (V2) -- +(270:5);
	\draw[thick,green!50!black, name path=V245] (V2) -- +(285:5);
	\path [name intersections={of=O15 and V115,by=A1}];
	\path [name intersections={of=O30 and V115,by=A2}];
	\path [name intersections={of=O45 and V115,by=A3}];
	\path [name intersections={of=O15 and V215,by=B1}];
	\path [name intersections={of=O30 and V215,by=B2}];
	\path [name intersections={of=O45 and V215,by=B3}];
	\path [name intersections={of=O30 and V130,by=C1}];
	\path [name intersections={of=O15 and V145,by=C2}];
	\path [name intersections={of=O15 and V130,by=C3}];
	\path [name intersections={of=V215 and V115,by=D1}];
	\path [name intersections={of=V245 and V145,by=D2}];
	\path [name intersections={of=V245 and O45,by=D3}];
	\path [name intersections={of=O45 and V145,by=D4}];
	\draw[blue, thick] ([shift={(-105:0.3cm)}]A1.center) -- ([shift={(-105:-0.5cm)}]A1.center);
    \draw[green!50!black, thick] ([shift={(-75:0.3cm)}]D4.center) -- ([shift={(-75:-0.5cm)}]D4.center);
	\draw[orange,fill] (O) circle (2.5pt);
	\draw[orange,fill] (V1) circle (2.5pt);
	\draw[orange,fill] (V2) circle (2.5pt);
	\draw[orange,fill] (A1) circle (2.5pt);
	\draw[orange,fill] (A2) circle (2.5pt);
	\draw[orange,fill] (A3) circle (2.5pt);
	\draw[white,fill] (B1) circle (1.5pt);
	\draw[orange,fill] (B2) circle (2.5pt);
	\draw[orange,fill] (B3) circle (2.5pt);
	\draw[orange,fill] (C1) circle (2.5pt);
	\draw[orange,fill] (C2) circle (2.5pt);
	\draw[orange,fill] (C3) circle (2.5pt);
	\draw[white,fill] (D1) circle (1.5pt);
	\draw[white,fill] (D2) circle (1.5pt);
	\draw[white,fill] (D3) circle (1.5pt);
	\draw[orange,fill] (D4) circle (2.5pt);
\end{tikzpicture}

\caption{In this real figure of the dual Hesse arrangement  two  lines (blue and green)  are poorly represented as  not connected. The are no double intersections, just triple intersections of lines at orange points.}
\label{Dual}
\end{figure}

Any attempt to illustrate  the dual Hesse arrangement on the \emph{real} plane  has some deficiency: either some  lines shall be represented as curved or broken or  not-connected, or  some points of the arrangement will  be missing. Our illustration for it is Figure \ref{Dual}.

\section{Effect of quadratic maps $Q_i$ on    on the dual Hesse arrangement}

A quadratic Cremona map $Q: \plano  \dashrightarrow \plano  $ with three non-collinear indetermination points amounts to the  blow up of the three points  $p_1, p_2, p_3$ and the contraction of  the three lines  joining the points.  These are involutive Cremona maps, so we can identify the set of  indetermination points of $Q$  with the set of indeterminations  of $Q^{-1} = Q$. 

The strict transforms of lines  by the quadratic map $Q$ may be:  an irreducible conic, if the line does not pass by any indetermination point; a line, if the line  passes by exactly one indeterminantion point; or a point, if the lines joins two indetermination points and therefore is contracted  by $Q$.

Since each  set $\mathcal{P}_3 (i)$ in the decomposition    
\[\mathcal{P}_{12} =   \mathcal{P}_3 (1)  \cup \mathcal{P}_3 (\tau)   \cup \mathcal{P}_3  (\tau^2) \cup \mathcal{P}_ 3 (\infty)\]
is composed by three non-collinear points, then each one may be the indetermination set of a quadratic Cremona map; in homogenous coordinates:
\[(\star\star\star)\begin{cases}   Q_1(x:y:z) = (y^2 -x z : x^2 - y z: z^2- x y  ) \\ Q_\tau (x:y:z)= (\tau y^2 - x z, \tau x^2 - y z: z^2- \tau^2 x y    ) \\ Q_{\tau^2} (x:y:z)= (\tau^2 y^2 - x z, \tau^2 x^2 - yz : z^2- \tau x y) \\ Q_\infty(x:y:z)=  ( y z: x z: x y )   \end{cases}   \]

We assert that each $Q_i$ ($i= 1, \tau, \tau^2, \infty$) preserves the dual Hesse arrangement:

\begin{pro}\label{Qiarranjo}
The   strict transform by  each $Q_i$, $i =1,\tau,\tau^2,\infty$,   of the set of lines  $\mathcal{L}_9$ of the dual Hesse arrangement    is the same set of lines  $\mathcal{L}_9$.  moreover, the effects of each $Q_i$  on the sets of three points $\mathcal{P}_3 (i)$ are  the following:

\[\begin{cases}
Q_1 (\mathcal{P}_3 (1) )= \mathcal{P}_3 (1),\quad  Q_1 (\mathcal{P}_3 (\tau) )= \mathcal{P}_3 (\tau^2),\quad  Q_1 (\mathcal{P}_3 (\infty) )= \mathcal{P}_3 (\infty )\\
Q_\tau (\mathcal{P}_3 (1) )= \mathcal{P}_3 (\tau^2),\quad Q_\tau (\mathcal{P}_3 (\tau) )= \mathcal{P}_3 (\tau),\quad  Q_\tau (\mathcal{P}_3 (\infty) )= \mathcal{P}_3 (\infty )\\
Q_{\tau^2} (\mathcal{P}_3 (1) )= \mathcal{P}_3 (\tau),\quad  Q_{\tau^2} (\mathcal{P}_3 (\tau^2) )= \mathcal{P}_3 (\tau^2),\,\, Q_{\tau^2} (\mathcal{P}_3 (\infty) )= \mathcal{P}_3 (\infty )\\
Q_\infty (\mathcal{P}_3 (\tau) )= \mathcal{P}_3 (\tau^2),\quad Q_\infty (\mathcal{P}_3 (1) )= \mathcal{P}_3 (1),\quad  Q_\infty (\mathcal{P}_3 (\infty) )= \mathcal{P}_3 (\infty )
\end{cases}\]

\end{pro}

Remark: although  the effects of $Q_1$ and $Q_\infty$ coincide in the four \emph{sets}  $\mathcal{P}_3 (i)$, they do not  coincide point to point. 

\begin{figure}
\begin{tikzpicture}[scale=1]
	\coordinate (O) at (0,0);
	\draw[thick] (O)-- (15:1);	
	\draw[thick, blue] (O) -- (30:1);
	\draw[thick, red] (O) -- (45:1);
	\coordinate (V1) at (5,0);	
	\draw[thick, blue] (V1) -- +(135:1);
	\draw[thick, red] (V1) -- +(150:1);
	\draw[thick] (V1) -- +(165:1);
	\coordinate (V2) at (60:5cm);
	\draw[thick,blue] (V2) -- +(255:5);
	\draw[thick,red] (V2) -- +(270:5);
	\draw[thick] (V2) -- +(285:5);
	\draw[thick,dashed,gray] (O) node[below]{$i$}-- (V1)node[below]{$i$} -- (V2)node[above]{$i$}  -- cycle node[midway,left,xshift=-5pt]{$\Delta_i$};
\end{tikzpicture}
\caption{A set  of points  $\mathcal{P}_3 (i)$ simply denoted by $i$, a  triangle $\Delta_i$ in dotted lines, and portions  of lines of $\mathcal{L}_{9}$ in red, blue, black colors.}
\label{Deltai}
\end{figure}
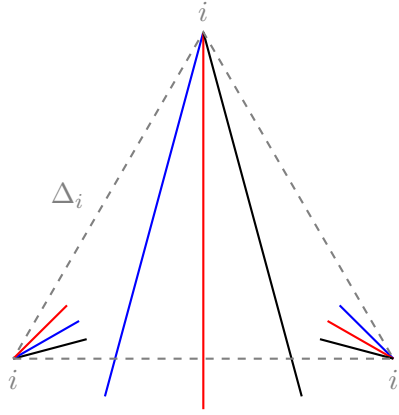

Proof: 

The effect of each quadratic map is  to  blow up the set $\mathcal{P}_3(i)$ and to contract  the lines of the triangle $\Delta_i$ connecting the points of $ \mathcal{P} (i)$: remark that no line of  such triangles  belongs to the dual Hesse arrangement.  Each line $l_{i, j}$ of $\Delta_i$  intersect six lines of $\mathcal{L}_{9}$ at two vertices of $\Delta_i$  and therefore intersects the remaing three lines of $\mathcal{L}_{9}$ out of these vertices.  Figure \ref{Deltai}  illustrates a system $\mathcal{P}_3(i)$ and  $\Delta_i$.
 
Each line of $\mathcal{L}_{9}$ pass by exactly one indetermination  point of the quadratic  $Q_i$, hence  its strict transform is again a line. And the contraction of each  $l_{i , j}$  of $\Delta_i$  brings together three strict transforms of lines of the arrangement.  

These facts imply that  the lines, their triple  point intersections  and the incidence relations $(12_3, 9_4  )$ of the arrangement are preserved.
The projective rigidity of the   dual Hesse arrangement assures that it  is preserved by each $Q_i$ (in the sense of struct transform).
 
In coordinates, for $\mathcal{L}_9:  (y^3-z^3) (y^3-x^3) (x^3-z^3)  = 0$ 
\[\begin{cases} Q_1^*( \mathcal{L}_9) =      (x+y+z)^3 (x+ \tau y + \tau^2 z  )^3  ( x + \tau^2 y + \tau z )^3 \cdot \mathcal{L}_9\\
Q_\tau^*( \mathcal{L}_9) =      (x + \tau y + \tau z)^3 (x+y + \tau^2 z)^3 (x+\tau^2 y + z)^3 \cdot  \mathcal{L}_9\\
 Q_{\tau^2}^*( \mathcal{L}_9) =      (x + \tau y + z)^3 (x + \tau^2 y + \tau^2 z)^3 (x+ y + \tau z)^3 \cdot  \mathcal{L}_9
 \end{cases}\]

Each set $\mathcal{P}_3 (i)$ is invariant by each involutive  $Q_i$. And the  effect on the other sets of three points  can be confirmed by means  of the effect of each $Q_i$ on lines  passing by  each point of the set, as  described  in ($\star\star$). 
\qed


\begin{figure}
   \begin{tikzpicture}[scale=1.5]
		\coordinate (O) at (0,0);
		\draw[orange,thick, name path=O15] (O) -- (15:5);	
		\draw[thick, name path=O30] (O) -- (30:5);
		\draw[orange,thick, name path=O45] (O) -- (45:5);
		\coordinate (V1) at (5,0);	
		\draw[green!50!black,thick, name path=V115] (V1) -- +(135:5);
		\draw[thick, name path=V130] (V1) -- +(150:5);
		\draw[green!50!black,thick, name path=V145] (V1) -- +(165:5);
		\coordinate (V2) at (canvas polar cs:radius=5cm,angle=60);
		\draw[thick,blue, name path=V215] (V2) -- +(255:5);
		\draw[thick, name path=V230] (V2) -- +(270:5);
		\draw[thick,blue, name path=V245] (V2) -- +(285:5);
		\path [name intersections={of=O15 and V115,by=A1}];
		\path [name intersections={of=O30 and V115,by=A2}];
		\path [name intersections={of=O45 and V115,by=A3}];
		\path [name intersections={of=O15 and V215,by=B1}];
		\path [name intersections={of=O30 and V215,by=B2}];
		\path [name intersections={of=O45 and V215,by=B3}];
		\path [name intersections={of=O30 and V130,by=C1}];
		\path [name intersections={of=O15 and V145,by=C2}];
		\path [name intersections={of=O15 and V130,by=C3}];
		\path [name intersections={of=V215 and V115,by=D1}];
		\path [name intersections={of=V245 and V145,by=D2}];
		\path [name intersections={of=V245 and O45,by=D3}];
		\path [name intersections={of=O45 and V145,by=D4}];
	\draw[blue, thick] ([shift={(-105:0.3cm)}]A1.center) -- ([shift={(-105:-0.5cm)}]A1.center);
\draw[blue, thick] ([shift={(-75:0.3cm)}]D4.center) -- ([shift={(-75:-0.5cm)}]D4.center);

		\draw[fill] (O) node[yshift=-2.5pt,below]{$\infty$} circle (1.5pt);
		\draw[fill] (V1) node[yshift=-2.5pt,below]{$\infty$} circle (1.5pt);
		\draw[fill] (V2) node[yshift=2.5pt,above]{$\infty$} circle (1.5pt);
		\draw[fill] (A1) node[xshift=2.5pt,yshift=-2.5pt,right]{$1$} circle (1.5pt);
		\draw[fill] (A2) node[right]{$\tau^2$} circle (1.5pt);
		\draw[fill] (A3) node[xshift=2.5pt,yshift=-2.5pt,right]{$\tau$} circle (1.5pt);
		\draw[white,fill] (B1) circle (1.5pt);
		\draw[fill] (B2) node[xshift=3pt,below]{$\tau$} circle (1.5pt);
		\draw[fill] (B3) node[yshift=-3pt,left]{$\tau^2$}  circle (1.5pt);
		\draw[fill] (C1) node[xshift=2.5pt,yshift=-2.5pt,below left]{$1$}  circle (1.5pt);
		\draw[fill] (C2) node[xshift=-6pt,below]{$\tau^2$}  circle (1.5pt);
		\draw[fill] (C3) node[xshift=-2pt,left]{$\tau$}  circle (1.5pt);
		\draw[white,fill] (D1) circle (1.5pt);
		\draw[white,fill] (D2) circle (1.5pt);
		\draw[white,fill] (D3) circle (1.5pt);
		\draw[fill] (D4) node[above,xshift=1pt]{$1$}  circle (1.5pt);
	\end{tikzpicture}

\caption{Effect of $Q_\infty$ on lines and points of the dual Hesse arrangement.  Points of $\mathcal{P}_3 (i)$ are denoted $i$. Black lines are invariant; at the center the fixed point $(1:1:1)$. Lines of the same color are switched by $Q_\infty$ }
\label{EfeitoQinfty}
\end{figure}
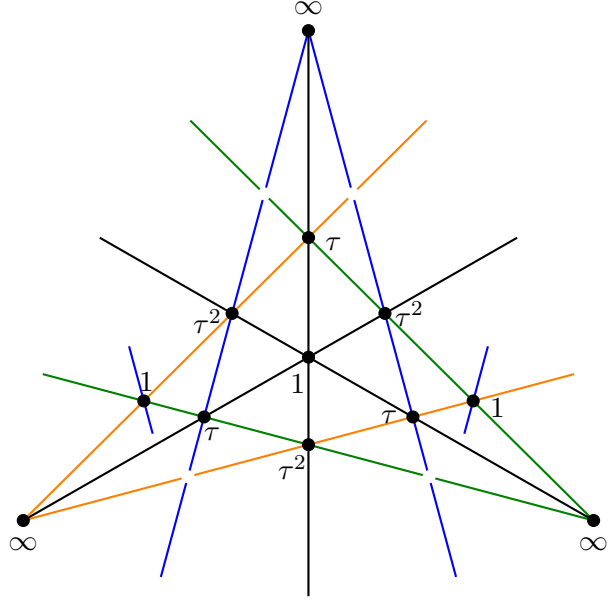

\begin{ex}

Let us examplify in the case $Q_\infty$ the fact that
\[ Q_\infty (\mathcal{P}_3 (\tau) )= \mathcal{P}_3 (\tau^2)  \]
(and reciprocally, since it is an involution).

As $Q_\infty (x:y:z)  = ( y z : x z : x y)$, in affine cordinates $(x,y)$ it is  given by $Q_\infty (x,y)= (\frac{1}{x}, \frac{1}{y})$ and we conclude that $(1,1)$ is a fixed point, that  the line $y=x$ is invariant, that the line $y = \tau x$ is sent to $y= \tau^2 x$, since $\tau^{-1} = \tau^2$, and vice-versa. That is, in the notation of $(\star)$, $l_ 1$ is invariant by  $Q_\infty$, $l_2$ is sent to $l_3$ (and reciprocaly). 

Examing in the other two affine coordinates, we conclude that $m_1$ is invariant, $m_2$ goes to $m_3$; that $n_1$ is invariant, $n_2$ goes to $n_3$. 

The effect of $Q_ \infty$ is  described  in  Figure \ref{EfeitoQinfty}.

\end{ex}

In general, the strict transform by a quadratic Cremona map  $Q$ with indeterminations $p_1,p_2,p_3$  of a curve $C$ of degree $\mbox{deg}(C)$  is a curve $\overline{C}$ of degree \[\mbox{deg}(\overline{C}) = 2 \, \mbox{deg}(C) - m(p_{j}, C) - m(p_{k},C) - m(p_{l},C)\]
  with points of multiplicites \[m(p_{j}, \overline{C}) = \mbox{deg}(C) - m_{p_k}(C)  - m_{p_l} (C), \quad j\neq k \neq l\]   

In the case of $Q_{i}$, for $i= 1,\tau,\tau^2,\infty$,  we can encode  its effect on curves as follows (recall that, thanks to Proposition \ref{Qiarranjo},  there are shifts in positions)

\begin{pro}\label{Qisimbolico} 

Let  $C$ be a curve of degree $d$,  having  the same multiplicity $m_1, m_\tau, m_{\tau^2}, m_{\infty} $ at the three points of the  $\mathcal{P}_3 (1), \mathcal{P}_ 3 (\tau), \mathcal{P}_ 3 (\tau^2), \mathcal{P}_3 (\infty)$,  respectively. These data represented by the list  \[[ d, m_1 , m_\tau, m_{\tau^2}, m_\infty]\]
Using the same convention,  the strict  transform $\overline{C}_i$  of  $C$ by each quadratic map $Q_i$  has degree and multiplicities:
\[\begin{cases} [ 2 d - 3 m_1,  d - 2  m_1,  m_{\tau^2}, m_ \tau , m_\infty],\quad   \mbox{if} \, \,  i= 1\\
[ 2 d - 3 m_\tau,   m_{\tau^2}, d - 2  m_\tau, m_1 , m_\infty],\quad   \mbox{if} \,\,  i= \tau\\
[ 2 d - 3 m_{\tau^2},   m_\tau ,  m_1 , d - 2  m_{\tau^2}, m_{\infty}] ,\quad   \mbox{if}\,    i = \,\,\,   \tau^2\\
[ 2 d - 3 m_\infty,   m_1,   m_{\tau^2}, m_\tau , d -  2   m_{\infty}],\quad   \mbox{if} \,\, \,    i = \infty\end{cases}\]

\end{pro}

\section{Lins Neto's  foliations $\F_t$ transformed by $Q_i$}\label{Ft}

Now we use the background concepts of Section \ref{background}, for instance, the   notions of \emph{degree} of a foliation $\mbox{deg}(\F)$ and the different  multiplicites $\mu(p, \F)$ and $l(p,\F)$.

The degree  as a divisor of the \emph{tangency set}  of a pair of degree  $d$ foliations in the plane  is $2 d +1$. This can be checked by contracting  polynomial   $1$-form and  polynomial vector field representing the pair of foliations.

In particular,  the degree of the tangency set of a pair of degree $4$ foliations is $9$. For this reason the  set of lines $\mathcal{L}_9$ of the dual Hesse arrangement can be   invariant by a pair of degre four foliations and by the   $1$-dimensional linear family of  foliations generated by the pair. This is the case, and  Lins Neto's family  of foliations
\[\F_t:\quad   \Omega + t \cdot  \Xi = 0, \quad  \mbox{for}   \quad  t \in \overline{\mathbb{C}} \quad \mbox{and} \]
{\small { \[\begin{cases} 
 \Omega  :=  z (y-z) (z^2+y^2+yz) y dx-z(x-z)(x^2+zx+z^2)x dy+xy(x-y)(x^2+xy+y^2)dz
\\ 
\Xi  := -(y-z)(z^2+y^2+yz)x^2   dx+(x-z)(x^2+zx+z^2)y^2dy-z^2(x-y)(x^2+xy+y^2) dz    
\end{cases} \]}}
is tangent to $\mathcal{L}_9$. 

In \cite{LinsNeto},  it is proved  that, except for $\F_{\infty}, \F_1,  \F_\tau, \F_{\tau^2}$, the singularities of $\F_t$ are $21$ distinct points:  twelve fixed radial points  at  the set $\mathcal{P}_{12}$; and  nine movable reduced nondegenerate points $p_1(\F_t),\ldots ,p_9(\F_t)$  with local first integrals of type $x^3 y =c$, one at each line of $\mathcal{L}_9$ out of the four points of $\mathcal{P}_{12}$ over the line.

 The effect on a  foliation $\F$ of $\plano$ of a quadratic  Cremona   map $Q$ with three  non-collinear indeterminations points $p_j, p_k, p_l$ is the following (cf.  \cite{MePe},  Lemma 1).
The  strict transformed foliation $\overline{\F}$ of $\F$ under $Q$ (i.e. with finite singular set)  has degree
\begin{equation}\label{changedegree} 
\mbox{deg}(\overline{\F})= 2 \cdot \mbox{deg}(\F)  + 2 -l(p_j,\F) - l(p_k,\F) -l(p_l,\F) 
\end{equation}
and multiplicities 
\begin{equation}\label{changelp}
l(p_k,\overline{\F}) = \mbox{deg}(\F)  + 2 - l(p_j,\F) - l(p_l,\F)
\end{equation}
For all $t \in \overline{\mathbb{C}}$ the non-reduced  singular points of the family $\F_t$ are either radial points  or    points  of local form $d( x y  (x-y) ) =0$:  in both cases,  $l(p,\F)= 2$. 

 Since the degree of $\F_t$ is four, we conclude from (\ref{changedegree}) that the strict transform of any $\F_t$ by $Q_i$ is again a degree  four foliation. And (\ref{changelp}) implies that the new singularities also have $l(p,\overline{F}) =2$. 

\begin{figure}
\begin{tikzpicture}[scale=1]
	\coordinate (O) at (0,0);
	\draw[thick] (O)-- (15:1);	
	\draw[thick, blue] (O) -- (30:1);
	\path[name path=O30] (O) -- (30:5);
	\draw[thick, red] (O) -- (45:1);
	\coordinate (V1) at (5,0);	
	\draw[thick, blue] (V1) -- +(135:1);
	\draw[thick, red] (V1) -- +(150:1);
	\path[name path=V130] (V1) -- +(150:5);
	\draw[thick] (V1) -- +(165:1);
	\coordinate (V2) at (60:5cm);
	\draw[thick,blue] (V2) -- +(255:1);
	\draw[thick,red] (V2) -- +(270:1);
	\draw[thick] (V2) -- +(285:1);
	\draw[thick,dashed,gray] (O) node[below]{$i$}-- (V1)node[below]{$i$} -- (V2)node[above]{$i$}  -- cycle node[midway,left,xshift=-5pt]{$\Delta_i$};
	\path [name intersections={of=O30 and V130,by=C1}];
	\node[green!50!black] at (C1) {$\mathcal{F}_t$};
	\foreach \t in {0.15,0.2,...,0.9}{
	   \coordinate (tmp) at ($(O)!\t!(V1)$);
	   \draw[thick,green!50!black] ([yshift=-0.2cm] tmp) to[bend right] ([yshift=0.2cm] tmp);
    }
	\foreach \t in {0.15,0.2,...,0.9}{
	\coordinate (tmp) at ($(V1)!\t!(V2)$);
	\draw[thick,green!50!black] ([xshift=-0.15cm,yshift=-0.15cm] tmp) to[bend right] (tmp) to[bend left] ([xshift=0.15cm,yshift=0.15cm] tmp);
}
	\foreach \t in {0.15,0.2,...,0.9}{
	\coordinate (tmp) at ($(V2)!\t!(O)$);
	\draw[thick,green!50!black] ([xshift=-0.15cm,yshift=0.15cm] tmp) to[bend left] (tmp) to[bend right] ([xshift=0.15cm,yshift=-0.15cm] tmp);
}
\foreach \pt in {O,V1,V2}{
\draw[thick,green!50!black] ([xshift=-0.3cm,yshift=-0.3cm] \pt) -- ([xshift=0.3cm,yshift=0.3cm] \pt);
\draw[thick,green!50!black] ([xshift=0.3cm,yshift=-0.3cm] \pt) -- ([xshift=-0.3cm,yshift=0.3cm] \pt);
\draw[thick,green!50!black] ([xshift=-0.1cm,yshift=-0.3cm] \pt) -- ([xshift=0.1cm,yshift=0.3cm] \pt);
}
\end{tikzpicture}
\caption{In green the foliation $\F_t$ and three radial points at vertices of $\Delta_i$. In red, blue, black local branches of $\mathcal{L}_{9}$ on the radial points}
\label{Deltaifol}
\end{figure}
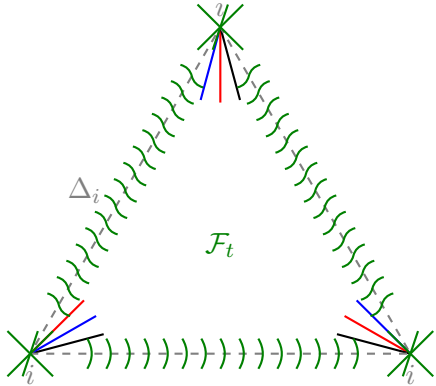

For $\F_t$ with  $t \notin \{i, \tau,\tau^2,\infty\}$,  the non-reduced singularities are   twelve radial points. Each line $l_{i, j}$ of the  fundamental triangle $\Delta_i$   of $Q_i$ passes  by two radial points. The sum of order of tangencies  of $l_{i, j}$ with $\F_t$ concentrated at the radial points  is $2 + 2 = 4$. Therefore, out of the radial points there is complete  transversality between   the foliations $\F_t$ and the lines $l_{i, j}$.  See  Figure \ref{Deltaifol}.

This implies that the points which arise as contraction of $l_{i, j}$ are again  radial points of the transformed foliation $\overline{\F_t}$.

As already remarked, the set of lines $\mathcal{L}_9$ of the dual Hesse arrangement  is preserved by  each $Q_i$. 

The conclusion is that, at least for  $t \notin \{1,\tau,\tau^2,\infty\}$,
$\F_t$ are transformed by $Q_i$  in another element $\F_{t^{\prime}}$ of the family.  We treat cases $t \in \{ 1,\tau,\tau^2,\infty\}$ separately in Example \ref{cubicos} and Section \ref{taut}.

\begin{pro}\label{efeitonoparametro}

 The  strict transform of   each foliation  of Lins Neto's family  $\F_t: \quad    \Omega + t \cdot \Xi  $ 
by each quadratic map  $Q_i$, $i \in \{1,\tau,\tau^2, \infty\}$ is another element $\F_{q_i(t)}:\quad    \Omega +  q_i(t)  \cdot \Xi $
where
\[\begin{cases} q_1  (t) = - t -1 \\
q_\tau(t)= -t - \tau\\
q_{\tau^2}(t) = -t - \tau^2\\
q_{\infty}(t) = \frac{1}{t}  \end{cases}\] 
 
\end{pro}

Proof:

It is a matter of computing  the pullbacks by $Q_i$ of the $1$-forms defining $\F_t$, extract common factors and get  from this the strict transforms $Q_ i(\F_t)$. 
In fact, 
\[Q_1^*( \Omega + t \cdot \Xi) = c\cdot (x+y+z)^2 (x+ \tau y + \tau^2 z  )^2  ( x + \tau^2 y + \tau z )^2 \cdot  (\Omega + (- t -1 ) \cdot  \Xi)\]
and therefore $q_1(t)= - t- 1$. 
\[Q_{\tau}^* (\Omega + t \cdot  \Xi) = c \cdot (x + \tau y + \tau z)^2(x+y + \tau^2 z)^2 (x+\tau^2 y + z)^2 \cdot (\Omega + (-t- \tau) \cdot   \Xi )\]
and threrefore  $q_\tau(t)= - t- \tau$.
\[Q_{\tau^2}^* (\Omega + t \cdot \Xi) = c \cdot (x + \tau y + z)^2(x + \tau^2 y + \tau^2 z)^2 (x+ y + \tau z)^2 \cdot (\Omega + (-t- \tau^2) \cdot  \Xi )\]
and threrefore  $q_{\tau^2} (t)= - t- \tau^2$. At last, 
\[Q_\infty^*(\Omega + t \cdot \Xi) = c\cdot  x^2 y^2 z^2 \cdot (t \cdot \Omega +  \Xi)\]
which means $q_\infty(t) = \frac{1}{t}$.

\qed

\section{Elliptic pencils resulting of applications of $Q_i$}\label{ellipticpencilssurfaces}

\begin{exe}\label{cubicos}
The four foliations  $\F_1,  \F_\tau, \F_{\tau^2}, \F_\infty$ are  elliptic    pencils of cubics. 

Each pencil $\F_i$   has   three special elements, which are   three  lines among those of  $\mathcal{L}_9$ concurrent at  the three points of $\mathcal{P}_3 (i)$.  

Each pencil $\F_i$ has  nine base-points, located at the  points $\mathcal{P}_{12}\setminus \mathcal{P}_3(i)$. 
As foliations, their   singular are   composed by  nine  radial points (at base-points of the pencil) and  three   points of local type $d( x \, y \, (x-y))=0$ (whose Milnor number is $4$, cf.   \cite{BrBGF} p. 5).

Abusing notation for the line and its equation, we can write the generators of each pencil, in the notation of ($\star$):
\[(\square)\begin{cases}
\F_1:\,   c_1 \cdot  l_1  m_1  n_1  + c_2 \cdot  l_2  m_3 n_2 = 0\\   
\F_\tau:\, c_1 \cdot  l_2 m_1 n_3 + c_2 \cdot  l_1 m_2 n_2 = 0 \\  
\F_{\tau^2}:\, c_1\cdot  l_1 m_3  n_3  + c_2 \cdot  l_2 m_2  n_1 =0\\ 
\F_{\infty}: c_1 \cdot  l_1  l_2  l_3  + c_2 \cdot m_1  m_2  m_3 =0 \end{cases}\]   

For  $\F_ 1$, the third special element  is $l_3 m_2 n_3 =0 $; for $\F _\tau$,  $l_3  m_3  n_1  =0$; for  $\F _{\tau^2}$, $ l_3  m_1  n_2 =0 $; for $\F_ \infty$, $ n_1 n_2  n_3  =0$.


Figure \ref{novof3g3Effective} illustrates schematicaly the pencil  $\F_\infty$ with triple points at the set  $\mathcal{P}_3 (\infty)$

As remarked in 
\cite{LinsNeto} (Prop. 2-c)   the pencils of cubics $\F_1, \F_\tau, \F_{\tau^2}, \F_\infty$  are  projectively equivalent.   We shall  return to the    projective equivalences of $\F_1, \F_\tau, \F_{\tau^2}$ and o the   projective equivalence of   $\F\infty$ an  $\F_1$  in Section \ref{taut}.

\begin{figure}
   \begin{tikzpicture}[scale=1.5]
	\coordinate (O) at (0,0);
	\draw[red,thick, name path=O15] (O) -- (15:5);	
	\draw[red,thick, name path=O30] (O) -- (30:5);
	\draw[red,thick, name path=O45] (O) -- (45:5);
	\coordinate (V1) at (5,0);	
	\draw[green!50!black,thick, name path=V115] (V1) -- +(135:5);
	\draw[green!50!black,thick, name path=V130] (V1) -- +(150:5);
	\draw[green!50!black,thick, name path=V145] (V1) -- +(165:5);
	\coordinate (V2) at (canvas polar cs:radius=5cm,angle=60);
	\draw[thick,blue, name path=V215] (V2) -- +(255:5);
	\draw[thick, blue, name path=V230] (V2) -- +(270:5);
	\draw[thick,blue, name path=V245] (V2) -- +(285:5);
	\path [name intersections={of=O15 and V115,by=A1}];
	\path [name intersections={of=O30 and V115,by=A2}];
	\path [name intersections={of=O45 and V115,by=A3}];
	\path [name intersections={of=O15 and V215,by=B1}];
	\path [name intersections={of=O30 and V215,by=B2}];
	\path [name intersections={of=O45 and V215,by=B3}];
	\path [name intersections={of=O30 and V130,by=C1}];
	\path [name intersections={of=O15 and V145,by=C2}];
	\path [name intersections={of=O15 and V130,by=C3}];
	\path [name intersections={of=V215 and V115,by=D1}];
	\path [name intersections={of=V245 and V145,by=D2}];
	\path [name intersections={of=V245 and O45,by=D3}];
	\path [name intersections={of=O45 and V145,by=D4}];
	\draw[blue, thick] ([shift={(-105:0.3cm)}]A1.center) -- ([shift={(-105:-0.5cm)}]A1.center);
\draw[blue, thick] ([shift={(-75:0.3cm)}]D4.center) -- ([shift={(-75:-0.5cm)}]D4.center);
	\draw[orange,fill] (O) node[black,yshift=-2.5pt,below]{$\infty$} circle (1.5pt);
	\draw[orange,fill] (V1) node[black,yshift=-2.5pt,below]{$\infty$} circle (1.5pt);
	\draw[orange,fill] (V2) node[black,yshift=2.5pt,above]{$\infty$} circle (1.5pt);
	\draw[orange,fill] (A1) node[black,xshift=2.5pt,yshift=-2.5pt,right]{$1$} circle (1.5pt);
	\draw[orange,fill] (A2) node[black,right]{$\tau^2$} circle (1.5pt);
	\draw[orange,fill] (A3) node[black,xshift=2.5pt,yshift=-2.5pt,right]{$\tau$} circle (1.5pt);
	\draw[white,fill] (B1) circle (1.5pt);
	\draw[orange,fill] (B2) node[black,xshift=3pt,below]{$\tau$} circle (1.5pt);
	\draw[orange,fill] (B3) node[black,yshift=-3pt,left]{$\tau^2$}  circle (1.5pt);
	\draw[orange,fill] (C1) node[black,xshift=2.5pt,yshift=-2.5pt,below left]{$1$}  circle (1.5pt);
	\draw[orange,fill] (C2) node[black,xshift=-6pt,below]{$\tau^2$}  circle (1.5pt);
	\draw[orange,fill] (C3) node[black,xshift=-2pt,left]{$\tau$}  circle (1.5pt);
	\draw[white,fill] (D1) circle (1.5pt);
	\draw[white,fill] (D2) circle (1.5pt);
	\draw[white,fill] (D3) circle (1.5pt);
	\draw[orange,fill] (D4) node[black,above,xshift=1pt]{$1$}  circle (1.5pt);
\end{tikzpicture}
\caption{Special elements of the pencil $\F_\infty$ represented in blue, red and green. The three  points of the sets $\mathcal{P}_3(i)$ are denoted just by  the symbols $i$.}
\label{novof3g3Effective}
\end{figure}
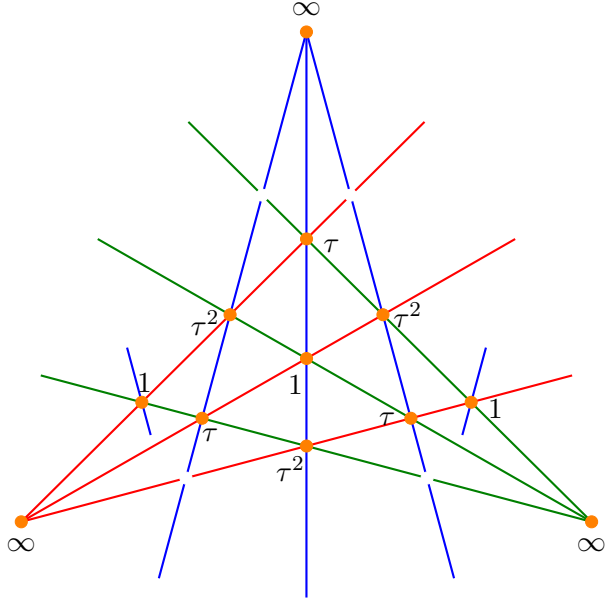

Now we  remark that   the pencils  $\F_1, \F_\tau, \F_{\tau^2}$  are strict transforms one of another by the quadratic Cremona maps $Q_i$. 
In fact, recalling that $\tau^2 = -1- \tau$ and $\frac{1}{\tau} = \tau^2$, it holds  (in the sense of strict transforms):
 \[\F_{\tau^2} =  Q_\tau (\F_1),\quad   \F_\tau = Q_\infty (\F_{\tau^2}),\quad   \F_\tau = Q_{\tau^2} (\F_1)\]
These three facts can be checked directly, by  applying  the quadratic transformations $Q_i$ listed in ($\star\star\star$)  to the pencils listed in ($\square$). 

Also we  remark that $Q_\infty (\F_1)= \F_1$ (in the sense of strict transform) can be checked directly in Figure \ref{EfeitoQinfty}: the three black lines of the arrangement through $(1:1:1)$ are invariant and the three lines through $(1:\tau:\tau^2)$ are switched with the three lines through $(1:\tau^2:\tau)$. Two of these three sets of lines generate $\F_ 1$ (recall $\P_3(1)= \{ (1:1:1),  (1:\tau:\tau^2), (1:\tau^2:\tau) \}$). 
 
\end{exe}

\begin{exe}\label{sextics}
The strict transform of $\F_\infty$ by $Q_\infty$ is $\F_0$, which is  an   elliptic pencil of   sextics . In fact
\[\F_ \infty: c_1 \cdot  (y-x) (y-\tau x) (y-\tau^2 x)+ c_2 \cdot (y-z)(y-\tau z) (y- \tau^2 z) =\]
\[=  c_1 \cdot (y^3-x^3) + c_2 \cdot (y^3-z^3) = 0\]
is sent by $Q_\infty (x:y:z)= (y z: x z: xy)$ to
\[\F_0:\,  c_ 1 \cdot ((xz)^3-(yz)^3) + c_2 \cdot ((xz)^3-(xy)^3) =\]
\[=  c_1 \cdot z^3 (x^3-y^3) + c_2 \cdot x^3 (z^3- y^3) = 0\]
There is a third special element in this pencil, $ (x^3-z^3) y^3  = 0$. See Figure \ref{novof6g6Effective}.

This pencil of sextics has  twelve base-points at the points  $\mathcal{P}_{12}$; its generic element is smooth at $\mathcal{P}_{12} \setminus \mathcal{P}_ 3(\infty)$ and has     ordinary $3$-uple points at $\{(0:0:1), (1:0:0), (0:1:0)\} = \mathcal{P}_ 3 (\infty) $. 

Remark  the agreement with Darboux' formula for pencils  in the plane  (cf. Section \ref{background}: $\F_0$ has three triple lines $x^3=0, y^3=0 , z^3 =0$ (i.e. $\alpha_{s,k} =3$) among its  special  elements $C_s$, therefore 
\[ 4  = \mbox{deg}(\F_0) = 2 \cdot 6 - 2 -  \underbrace{3 \cdot 1 \cdot (3-1)}_{\mbox{3 triple  lines}}\]



\begin{figure}
   \begin{tikzpicture}[scale=1.5]
	\coordinate (O) at (0,0);
	\draw[red,thick, name path=O15] (O) -- (15:5);	
	\draw[red,thick, name path=O30] (O) -- (30:5);
	\draw[red,thick, name path=O45] (O) -- (45:5);
	\coordinate (V1) at (5,0);	
	\draw[green!50!black,thick, name path=V115] (V1) -- +(135:5);
	\draw[green!50!black,thick, name path=V130] (V1) -- +(150:5);
	\draw[green!50!black,thick, name path=V145] (V1) -- +(165:5);
	\coordinate (V2) at (canvas polar cs:radius=5cm,angle=60);
	\draw[thick,blue, name path=V215] (V2) -- +(255:5);
	\draw[thick, blue, name path=V230] (V2) -- +(270:5);
	\draw[thick,blue, name path=V245] (V2) -- +(285:5);
	\draw[thick,green!50!black] (O) -- (V2) node[yshift=-2pt,midway,left]{$(z^3=0)$};
	\draw[thick,blue] (O) -- (V1) node[midway,below left]{$(y^3=0)$};
	\draw[thick,red] (V1) -- (V2) node[midway,right,yshift=-2pt]{$(x^3=0)$};
	\path [name intersections={of=O15 and V115,by=A1}];
	\path [name intersections={of=O30 and V115,by=A2}];
	\path [name intersections={of=O45 and V115,by=A3}];
	\path [name intersections={of=O15 and V215,by=B1}];
	\path [name intersections={of=O30 and V215,by=B2}];
	\path [name intersections={of=O45 and V215,by=B3}];
	\path [name intersections={of=O30 and V130,by=C1}];
	\path [name intersections={of=O15 and V145,by=C2}];
	\path [name intersections={of=O15 and V130,by=C3}];
	\path [name intersections={of=V215 and V115,by=D1}];
	\path [name intersections={of=V245 and V145,by=D2}];
	\path [name intersections={of=V245 and O45,by=D3}];
	\path [name intersections={of=O45 and V145,by=D4}];
	\draw[blue, thick] ([shift={(-105:0.3cm)}]A1.center) -- ([shift={(-105:-0.4cm)}]A1.center);
\draw[blue, thick] ([shift={(-75:0.3cm)}]D4.center) -- ([shift={(-75:-0.4cm)}]D4.center);
	\draw[fill] (O) node[black,yshift=-2.5pt,below]{$\infty$} circle (1.5pt);
	\draw[fill] (V1) node[black,yshift=-2.5pt,below]{$\infty$} circle (1.5pt);
	\draw[fill] (V2) node[black,yshift=2.5pt,above]{$\infty$} circle (1.5pt);
	\draw[fill] (A1) circle (1.5pt);
	\draw[fill] (A2) circle (1.5pt);
	\draw[fill] (A3) circle (1.5pt);
	\draw[white,fill] (B1) circle (1.5pt);
	\draw[fill] (B2) circle (1.5pt);
	\draw[fill] (B3) circle (1.5pt);
	\draw[fill] (C1) circle (1.5pt);
	\draw[fill] (C2) circle (1.5pt);
	\draw[fill] (C3) circle (1.5pt);
	\draw[white,fill] (D1) circle (1.5pt);
	\draw[white,fill] (D2) circle (1.5pt);
	\draw[white,fill] (D3) circle (1.5pt);
	\draw[fill] (D4) circle (1.5pt);
\end{tikzpicture}
\caption{The special elements of the pencil of sextics $\F_0$ in red, blue, black.}
\label{novof6g6Effective}
\end{figure}
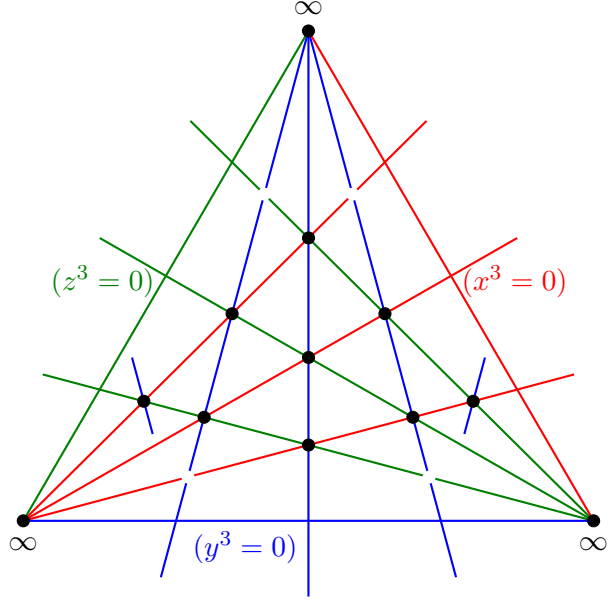

\end{exe}

\begin{exe}\label{nonics}

The strict transform of the pencil of sextics $\F_ 0$  (of previous Example \ref{sextics})  by $Q_1$ is a pencil of nonics  $\F_ {-1} = Q_1 (\F_ 0)$ 
\[\F_ {-1}:\, c_1\cdot (y^3-z^3)(x^2-y z)^3 + c_2 \cdot (x^3-z^3)(-y^2+zx)^3 = 0\] 
The generic element  has  $4$-uple ordinary points at 
\[ \{ (1:1:1), (1: \tau, \tau^2), (1:\tau^2,\tau)\} = \mathcal{P}_3 (1),  \]
$3$-uple points at
\[ \{ (0:0:1), (1:0:0), (0:1:0) \} = P_3 (\infty)\ \]
 and smooth points at  the extra six points of $\mathcal{P}_{12}$. 

There are three special elements; each one is composed by three lines and a triple conic (image by  $Q_1$ of a triple line of $\F_0$). Illustrated by Figure \ref{novof9g9Effective}.

Remark the agreement with Darboux' formula for pencils in the plane  ($\alpha_{s,k} =3$) 
\[ 4 = \mbox{deg}(\F_0) = 2 \cdot 9  - 2 - \underbrace{ 3 \cdot 2 \cdot  (3-1)}_{\mbox{3  triple conics}}\]

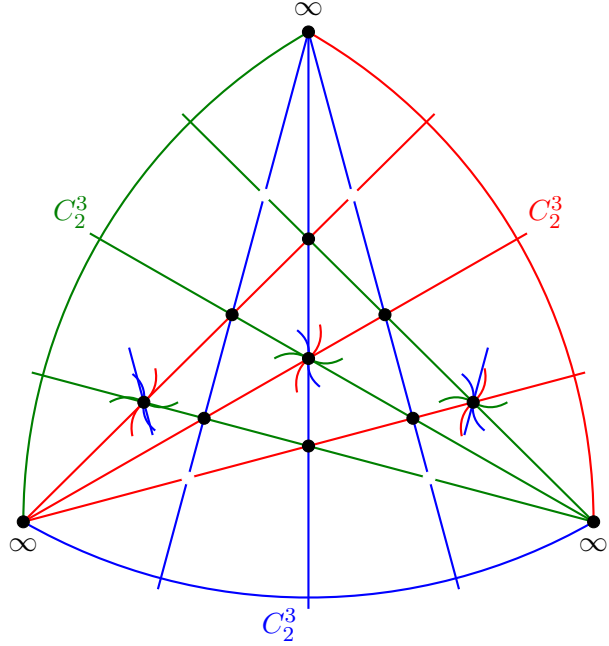
\begin{figure}
   \begin{tikzpicture}[scale=1.5]
	\coordinate (O) at (0,0);
	\draw[red,thick, name path=O15] (O) -- (15:5.1);	
	\draw[red,thick, name path=O30] (O) -- (30:5.1);
	\draw[red,thick, name path=O45] (O) -- (45:5.1);
	\coordinate (V1) at (5,0);	
	\draw[green!50!black,thick, name path=V115] (V1) -- +(135:5.1);
	\draw[green!50!black,thick, name path=V130] (V1) -- +(150:5.1);
	\draw[green!50!black,thick, name path=V145] (V1) -- +(165:5.1);
	\coordinate (V2) at (canvas polar cs:radius=5cm,angle=60);
	\draw[thick,blue, name path=V215] (V2) -- +(255:5.1);
	\draw[thick, blue, name path=V230] (V2) -- +(270:5.1);
	\draw[thick,blue, name path=V245] (V2) -- +(285:5.1);
\draw[thick,red] ([shift=(0:5cm)]O) arc (0:60:5cm) node[midway,above right]{$C_2^3$}; 
\draw[thick,green!50!black] ([shift=(120:5cm)]V1) arc (120:180:5cm) node[midway,above left]{$C_2^3$};
\draw[thick,blue] ([shift=(240:5cm)]V2) arc (240:300:5cm) node[midway,below left]{$C_2^3$};
	\path [name intersections={of=O15 and V115,by=A1}];
	\path [name intersections={of=O30 and V115,by=A2}];
	\path [name intersections={of=O45 and V115,by=A3}];
	\path [name intersections={of=O15 and V215,by=B1}];
	\path [name intersections={of=O30 and V215,by=B2}];
	\path [name intersections={of=O45 and V215,by=B3}];
	\path [name intersections={of=O30 and V130,by=C1}]; 
	\path [name intersections={of=O15 and V145,by=C2}];
	\path [name intersections={of=O15 and V130,by=C3}];
	\path [name intersections={of=V215 and V115,by=D1}];
	\path [name intersections={of=V245 and V145,by=D2}];
	\path [name intersections={of=V245 and O45,by=D3}];
	\path [name intersections={of=O45 and V145,by=D4}];
	\draw[blue, thick] ([shift={(-105:0.3cm)}]A1.center) -- ([shift={(-105:-0.5cm)}]A1.center);
\draw[blue, thick] ([shift={(-75:0.3cm)}]D4.center) -- ([shift={(-75:-0.5cm)}]D4.center);
\foreach \x in {C1,A1,D4} {
	\draw[blue,thick] ([xshift=0.1cm,yshift=-0.25cm] \x) to[bend left] (\x) to[bend right] 
	([xshift=-0.1cm,yshift=0.25cm] \x);
	\draw[red,thick] ([xshift=0.1cm,yshift=0.3cm] \x) to[bend left] (\x) to[bend right] ([xshift=-0.1cm,yshift=-0.3cm] \x);
	\draw[green!50!black,thick] ([xshift=0.3cm,yshift=0cm] \x) to[bend left] (\x) to[bend right] ([xshift=-0.3cm,yshift=0cm] \x);
}
	\draw[fill] (O) node[black,yshift=-2.5pt,below]{$\infty$} circle (1.5pt);
	\draw[fill] (V1) node[black,yshift=-2.5pt,below]{$\infty$} circle (1.5pt);
	\draw[fill] (V2) node[black,yshift=2.5pt,above]{$\infty$} circle (1.5pt);
	\draw[fill] (A1) circle (1.5pt);
	\draw[fill] (A2) circle (1.5pt);
	\draw[fill] (A3) circle (1.5pt);
	\draw[white,fill] (B1) circle (1.5pt);
	\draw[fill] (B2) circle (1.5pt);
	\draw[fill] (B3) circle (1.5pt);
	\draw[fill] (C1) circle (1.5pt);
	\draw[fill] (C2) circle (1.5pt);
	\draw[fill] (C3) circle (1.5pt);
	\draw[white,fill] (D1) circle (1.5pt);
	\draw[white,fill] (D2) circle (1.5pt);
	\draw[white,fill] (D3) circle (1.5pt);
	\draw[fill] (D4) circle (1.5pt);

\end{tikzpicture}
\caption{The special elements of the pencil of nonics $\F_{-1}$ in red, blue, green. Curved lines are representation of conics, poorly represented as not connected.}
\label{novof9g9Effective}
\end{figure}



\end{exe}

What  was remarked in Examples \ref{sextics} and \ref{nonics}  for $n=2,3$ continues to hold  under applications of $Q_i$. Are produced   elliptic pencils of curves of degrees $3 n$, having  three special elements, each one composed  by three lines and a rational curve of degree $n-1$ taken  with multiplicity three.

The agreement with Darboux' formula for pencil can be checked as:
\[ 4 =  2 \cdot 3 n - 2 -  \underbrace{3 \cdot (n-1) \cdot (3-2),}_{\mbox{3 triple curves of degree   n-1}} \quad n= 2, 3, \ldots\]

The rational curves which are supports of the  multiplicity three components  give rise to $(-1)$-curves  in the blown up plane in the twelve points $\mathcal{P}_{12}$ of the dual Hesse arrangement. These $(-1)$ curves are components of singular fibers of non-minimal elliptic fibrations. And this topic  shall be develloped at lenght in \cite{MPNegative}.

\section{Proof of Theorem \ref{principal}}

The effect  of the quadratic maps $Q_i$ on the parameter $t= m + n \tau$ (cf. Proposition \ref{efeitonoparametro}) of the familiy $\F_ t$ is 
\[\begin{cases} q_1  (t) = - t -1 \\
q_\tau(t)= -t - \tau\\
q_{\tau^2}(t) = -t - \tau^2 = -t + 1 +\tau \\ q_{\infty}(t) = \frac{1}{t}\end{cases}\]
For $t= m+ n \tau = (m,n) \in \mathbb{Z} \times \mathbb{Z}$,  we write
\[\begin{cases} q_1((m,n)) = (-m , -n ) + (-1,0)= (-m -1, -n) \\
q_\tau ((m,n)) = (-m , -n ) + (0,-1)= (-m, -n-1)\\
q_{\tau^2}((m,n)) = (-m , -n ) + (1,1)= (-m + 1, -n+1)  
\end{cases}\]

We introduce a norm in $\mathbb{Z}(\tau)$ given by 
\[N(m,n) := m^2+ n^2\]
and we have:

\begin{pro}\label{decreasingnorm}
\[ \begin{cases} N(q_1((m,n))) < N(m,n) \iff   m  \leq -1  \\ N(q_\tau ((m,n))) < N(m,n) \iff   n  \leq -1 \\ N(q_{\tau^2}((m,n))) < N(m,n) \iff   1 < m+n    \end{cases} \]

\end{pro}

Proof:

in fact, 
\[N(q_1((m,n))) < N(m,n) \iff  m^2 + n^2 + 2 m +1 < m^2 + n^2,\]
that is, $m\leq -1$ for  the integer $m$.
\[N(q_{\tau}((m,n))) < N(m,n) \iff  m^2 + n^2 + 2 n +1 < m^2 + n^2,\]
that is, $n\leq -1$ for the  integer $n$.
\[N(q_{\tau^2}((m,n))) < N(m,n) \iff  m^2 + n^2 - 2 m - 2 n +2 < m^2 + n^2,\]
that is, $ 1 < m +n $.

\qed 

Proposition \ref{decreasingnorm} is the basis for our algorithm:  it tells which $q_ 1, q_\tau, q_{\tau^2}$ to be applied in order to decrease the norm $N(m,n)$ of $t= m + n \tau$, $m , n \in \mathbb{Z}$.  

After a finite number of applications, we get  the conditon $N(m,n) \leq 1$.   

Next we give the list of nine pencils $\F_ {m + n \tau}$ for which $|m|\leq 1$  and $|n| \leq 1$. The generators os the pencils are presented in the form of  first integrals of the foliations:

\begin{pro}\label{fundamental}

 The list  $[d, m_1, m_\tau, m_{\tau^2}, m_\infty]$ give the degree $d$ of the generic element of the pencil and its  multiplicities of  at the three points of the sets $\mathcal{P}(1),   \mathcal{P}(\tau), \mathcal{P}(\tau^2), \mathcal{P}(\infty)$, respectively. 

\[\F_1 :   [3 , 0, 1,1,1],\,  \frac{(y-x)(z-x) (y-z)}{(y-\tau x)(z-\tau^2 x) (z-\tau y)} = c \]
\[\F_\tau :  [3 , 1, 0 ,  1, 1 ],\, \frac{ (y- \tau^2 x)(z-x) (z-\tau^2 y)}{ (y-x) (z-\tau x) (z-\tau y)} = c \]
\[\F_{\tau^2 }= \F _{-1-\tau}: [ 3, 1,1,0,1],\,   \frac{(y-x)(z-\tau^2 x) (z-\tau^2 y) }{(y-\tau x) (z-\tau x) (z-y) } = c \] 
\[\F_0:  [6,1,1,1,3], \quad   \frac{y^3 (x^3-z^3)}{x^3 (y^3-z^3)} = c\]
\[\F_{-1}: [9,4,1,1,3],\quad \frac{(y^3-z^3)(x^2-yz)^3}{(x^3-z^3)(-y^2+zx)^3} =c\]
\[\F_{-\tau}: [9,1,4,1,3],\quad  \frac{(y^3-z^3)(x^2- \tau^2 y z)^3}{(x^3-z^3)(-x z + \tau y^2)^3} = c \]
\[\F_{1+ \tau}:  [ 9, 1,1,4,3  ],  \frac{(y^3-z^3)(x^2- \tau y z)^3  }{(x^3-z^3) (y^2- \tau x z)^3 } =c \]
\[\F_{1- \tau}:  [15, 1, 7, 4 , 3    ],  \]
\[  \frac{(x^3-z^3) (x^3 y+\tau^2 y^4 + (-\tau+1) x y^2 z + (-\tau+1) x^2 z^2+y z^3)^3  }{(y^3- z^3) (x^4+\tau x y^3+ (2 \tau +1) x^2 y z+(2\tau+1)y^2 z^2+\tau x z^3)^3  } = c   \] 
\[\F_{-1 + \tau}:  [ 15, 7, 1,4,3  ],\,\]
\[ \frac{(x^3- z^3) (x^3 y+\tau y^4+ (-2 \tau -1) x y^2 z+(\tau-1) x^2 z^2+ y z^3 )^3 }{(y^3-z^3) (x^4+ (-\tau -1) x y^3+(\tau-1) x^2 y z+(\tau+2) y^2 z^2+(-\tau -1) x z^3  )^3 } = c\]

\end{pro}

Now we  remak that all nine pencils above can be  obtained  from  $\F_1$ and  $\F_{\infty}$ by applying $Q_i$.

In fact, in the strict transform sense,  it holds:
\[ \F_0 = Q_\infty(\F_\infty),\quad \F_{1+ \tau} = Q_{\tau^2}(\F_0),\]
\[\quad \F_{-1} = Q_1(\F_0),\quad  \F_{-\tau}= Q_\tau (\F_0),\]
\[  \F_ {-1+\tau} = Q_1(\F_ {-\tau}), \quad  \F_ {1-\tau} = Q_ \tau (\F_{-1}),\]
\[ \F_{-1-\tau}= \F_ {\tau^2} = Q_\tau (\F_1)\]
\[\F_\tau = Q_{\tau^2}(\F_1),\]
where the last assertion follows from \[q_ {\tau^2} (t)= -t- \tau^2= - t + 1 + \tau\]  applied to $t=1$.

Our algorithm transforms  the   pencils  $\F_ {m +n \tau}$, by repeated applications of $Q_ 1, Q_\tau, Q_ {\tau^2},  Q_ \infty$,   to either the  pencil of cubics   $\F_1$ or to  $ \F_\infty$. The data of these pencils  are
\[\F_1:\,  [3,0,1,1,1]   \]
(the generic smooth cubic of $\F_1$ does not pass by $\mathcal{P}_3(1)$), 
\[\F\infty:\, [3, 1,1,1,0]\] 
(the generic smooth cubic of $\F_\infty$ does not pass by $\mathcal{P}_3(\infty)$).
 
When arriving at one of theses two  pencils of cubics, it is a matter of inverting the order of the quadratic maps which were used, and we  obtain $\F_ {m + n \tau}$ step by step   from the pencils of cubics.



Remark: After having determined the data 
$[d, m_1,m_\tau,m_{\tau^2},m_{\infty}]$  of the elliptic pencil $\F_{t=m + n \tau}$, for  $m, n \in \mathbb{Z}$, the data of the pencil $\F_ {t= \frac{1}{m+n \tau}}$ can be obtained after one application of $q_\infty$, that is,
\[[2 d - 3 m_\infty, m_1, m_{\tau^2}, m_\tau, d - 2 m_\infty]   \]

\section{Examples and Algorithms in Python and Singular}\label{sectionexamples}

\begin{exe}

For $t=  - 2 + 8 \tau$, the algorith runs as follows:

\[\begin{cases} q_{\tau^2}: (-2,8) \mapsto ((2,-8) + (1,1) = (3,-7)\\
q_\tau:  ( 3,-7) \mapsto  (-3,7) + (0,-1) = (-3,6)\\
q_ {\tau^2}: (-3,6) \mapsto (3,-6) +(1,1) = (4,-5)\\
q_\tau:  ( 4,-5) \mapsto  (-4,5) + (0,-1) = (-4,4)\\
q_1:  (-4,4) \mapsto   (4,-4) + 9-1,0) = (3,-4)\\
q_\tau:  (3,-4) \mapsto  ((-3,4) + (0,-1) = (-3,3)\\
q_1: (-3,3) \mapsto (3,-3) + (-1,0) = (2,-3)\\
q_\tau: (2,-3) \mapsto (-2,3) + (0,-1) = (-2,2)  \\
q_1:  (2,-2)  \mapsto (2,-2)+ 9-1,0)= (1,-2)\\
q_\tau: (1,-2) \mapsto (-1,2) + (0,-1) = (-1,1)\\
q_1 :(-1,1) \mapsto (1,-1)+ (-1,0) = (0,-1)\\
q_\tau: (0,1) \mapsto (0,1) +(0,-1) = (0,0)\end{cases}\]
and finally  $\F_0 = Q_\infty (\F_\infty)$. 

The pencil $\F_\infty$ has data $[3,1,1,1,0]]$
applying to it the sequence of quadratic (in reversed order)  we obtain the sequence of pencils 
\[[6, 1, 1, 1, 3], [9, 1, 4, 1, 3],[15, 7, 1, 4, 3],
[27, 4, 13, 7, 3],\]
\[  [42, 19, 7, 13, 3], [63, 13, 28, 19, 3],[87, 37, 19, 28, 3],[117, 28, 49, 37, 3],\]
\[  [150, 61, 37, 49, 3], [189, 49, 76, 61, 3], [195, 76, 49, 67, 3], [243, 67, 97, 76, 3],\]
\[ [258, 97, 67, 91, 3]\]

This is the example default in the \emph{Python} algorithm   below, which presents the quadratic maps using 0 for $Q_\infty$,  $1$ for $Q_1$, $2$ for $Q_\tau$, $3$ for $Q_{\tau^2}$  and shows all intermediary pencils form the the cubic up to  to the pencil of  degree $258$. 
\end{exe}
\pagebreak

\begin{lstlisting}

m = -2; n = 8;
listmn=[m,n]
print('initial values  m,n of t=m+ntau:')
print(listmn)
listQ=[]
while  abs(listmn[0])+ abs(listmn[1])>1: 
    if listmn[0] + listmn[1] >1: 
        listmn[0] = -listmn[0] +1
        listmn[1]=  -listmn[1] +1
        listQ.append(3)
    elif  listmn[0]<=  listmn[1]:
        listmn[0]= -listmn[0] -1
        listmn[1]= -listmn[1]
        listQ.append(1)
    elif listmn[1]< listmn[0]:
        listmn[0]= -listmn[0]
        listmn[1]= -listmn[1]-1
        listQ.append(2)
else:   #  |m|  + |n|<= 1  reach the fundamental pencils
    if listmn[0]==0 and listmn[1]==0:
        listQ.append(0)
        fundamental=[ 3,1,1,1,0] #F\infty
    elif listmn[0]==-1 and listmn[1]==0:
        listQ.append(1)
        listQ.append(0)
        fundamental=[ 3,1,1,1,0] #F\infty
    elif listmn[0]==1 and listmn[1]==0:
        fundamental = [3,0,1,1,1]  #F_1
    elif listmn[0]==0 and listmn[1]==-1:
        listQ.append(2)
        listQ.append(0)
        fundamental =[3,1,1,1,0] #F_\infty
    elif listmn[0]==0 and listmn[1]==1:
        listQ.append(3)
        fundamental = [3,0,1,1,1] #F_1
def q1(list):  #q1
    l0=list[0];l1= list[1];l2=list[2];l3=list[3];
    list[0]=2*l0-3*l1
    list[1]= l0-2*l1
    list[2]=l3
    list[3]=l2
    return(list)
def q2(list): # q_tau
    l0=list[0];l1= list[1];l2=list[2];l3=list[3];
    list[0]= 2*l0-3*l2
    list[1]=l3
    list[3]=l1
    list[2]=l0-2*l2
    return(list)
def  q3(list): #q_tau^2
    l0=list[0];l1= list[1];l2=list[2];l3=list[3];
    list[0]= 2*l0-3*l3
    list[1]= l2
    list[2]=l1
    list[3]=l0-2*l3
    return(list)
def q0(list):  #q_infty
    l0=list[0];l2=list[2];l3=list[3];l4=list[4];
    list[0]= 2*l0-3*l4
    list[2]=l3
    list[3]=l2
    list[4]=l0-2*l4
    return(list)
list=fundamental #starting with the fundamental,  aplly Q's
print('and corresponding fundamental pencil:')
print(fundamental)
print('list of Qs to be applied to the fundamental pencil (from left to right),')
print('where 0 means Q_infty, 1 means Q_1, 2 means  Q_tau, 3 means  Q_tau^2:')
listQ.reverse()
print(listQ)
func_q_list = [q0,q1,q2,q3]
print('and effects of Qs: ')
for j in listQ:
     print(func_q_list[j](list))

\end{lstlisting}

\begin{exe}

For $t= -2 + 41 \tau$, acording to \cite{Puchuri}  $d = 5307$. According to our algoritm, starting with $\F_{\infty} =   [3,1,1,1,0] $ 
and applying  57  quadratic Cremona maps as in the list 
\[[0, 2, 1, 2, 1, 2, 1, 2, 1, 2, 1, 2, 1,\]
\[ 2, 1, 2, 1, 2, 1, 2, 1, 2, 1,\]
\[  2, 1, 2, 1, 2, 1, 2, 1, 2, 3, \]
\[ 2, 3, 2, 3, 2, 3, 2, 3, 2, 3, 2,\]
\[ 3, 2, 3, 2, 3, 2, 3, 2, 3, 2, 3, 2, 3]\]
(from lef to right, $0 = Q_\infty$, $1 = Q_1, 2= Q_\tau, 3= Q_{\tau^2}$), 
we get, step by step,
\[[6, 1, 1, 1, 3],[9, 1, 4, 1, 3], [15, 7, 1, 4, 3],[27,4, 13, 7, 3],[42, 19, 7, 13, 3],\]
\[[63, 13, 28, 19, 3,[87, 37, 19, 28, 3],[117, 28, 49, 37, 3],,[150, 61, 37, 49, 3],\]
\[[189, 49, 76, 61, 3],[231, 91, 61, 76, 3],[279, 76, 109, 91, 3],[330, 127, 91, 109, 3],\]
\[[387, 109, 148, 127, 3], [447, 169, 127, 148, 3],[513, 148, 193, 169, 3],[582, 217, 169, 193, 3]\]
\[[657, 193, 244, 217, 3],[735, 271, 217, 244, 3],[819, 244, 301, 271, 3]\]
\[[906, 331, 271, 301, 3],[999, 301, 364, 331, 3],[1095, 397, 331, 364, 3]\]
\centerline{    etc ... etc }
\[[4842, 1567, 1693, 1579, 3], [4947, 1693, 1567, 1684, 3],[5193, 1684, 1813, 1693, 3],\]
\[[5307, 1813, 1684, 1807, 3]\]

\end{exe}

For obtaining the expressions  of generators of the elliptic pencils in a  most simple form, we factorize the composition of the initial  cubic expressions with the $Q_i$ using the  \emph{Singular} software - which can be used on line  at $\mbox{https://www.singular.uni-kl.de/ }$.

For instance, consider 
\[\F_{-2-2\tau} =  Q_\tau (Q_{\tau^2} (Q_1 (Q_\infty(\F_\infty))))\]
which has data of degree and multiplciites $[18,7,7,1,3]$.

Next code in \emph{Singular} provide the irreducible factors (and a list of multiplicites of each factor) of the total transforms by any  sequence of quadratic maps. The default example start with   generators of $\F_ \infty$, written as $f=y^3- x^3$ and $g= y^3-z^3$, and apply the sequence defining $\F_{-2-2\tau}$. Remark that we used $a =\tau$ in the code.

After eliminating common factors (of degree one and two in this case) we get the generators of $\F_{-2-2\tau}$:
\[\F_{-2-2\tau}:\, c_1 \cdot  (x^3-y^3)  (f_5)^3 + c_2 \cdot  (x^3- z^3) (g_5)^3 = 0\]
where 
\[f_5 = x^4y+xy^4+3 \tau x^2y^2z+(-2\tau -2)x^3z^2+(-2\tau-2)y^3z^2+xyz^3+(\tau+1)z^5,\]
\[g_5= x^3y^2-\frac{1}{2}y^5+\frac{\tau}{2}x^4z+\frac{\tau}{2}xy^3z -\frac{3}{2}(\tau +1)x^2yz^2+y^2z^3+\frac{\tau}{2}xz^4  \]

\begin{lstlisting}

ring R=(0,a),(x,y,z),dp;                                      
minpoly=a2+a+1;                                                     
poly f = y^3-x^3;  
poly g = y^3-z^3;                                                                                   
map Qinf = R, y*z, x*z,  x*y;                                                                    
map Q1 = R, y^2- x*z, x^2- y*z, z^2- x*y;                                                        
map Qt = R, a*y^2- x*z, a*x^2- y*z, z^2- a^2 *x*y;                                               
map Qt2 = R, a^2*y^2- x*z, a^2 *x^2- y*z, z^2- a*x*y; 
map Cr= Qinf;
Cr=Q1(Cr);
Cr=Qt2(Cr);
Cr= Qt(Cr);
poly Crf=Cr(f); 
poly Crg =Cr(g); 
factorize(Crf);
factorize(Crg);

\end{lstlisting}  

\section{On the trivolution   $\tau \cdot t$ and   involutions  $ \frac{t+2}{t-1}$ and $-t$}\label{taut}

\begin{pro}

The linear automorphisms 
\[\begin{cases} T_1(x:y:z) = (z:x:y)\\   T_2(x:y:z) = (x: \tau y: \tau^2 z)\\   T_3(x:y:z) = (x: \tau^2 y:  \tau z)\end{cases}\]
 act as identity on the parameter of $\F_t$, that is,  $ t\mapsto t$. The linear automorphism of order three \[T_\tau (x:y:z)= (x: \tau y : z)\]
 produces in the plane
 \[T_\tau (\mathcal{P}_3 (1) ) = \mathcal{P}_3 (\tau), \quad T^2_\tau (\mathcal{P}_3 (1) ) = \mathcal{P}_3 (\tau^2)  \]
 acts on the parameter of $\F_t$ as $t \mapsto \tau^2 t$; and so $T\tau \circ T_\tau $ acts as $t \mapsto \tau \cdot t $ 

\end{pro}

Proof:

Just a matter of taking  pullbacks by these automorphisms of 
\[\Omega  + t \cdot E = 0\]
defining the foliations $\F_t$ (see Section \ref{Ft}).

\qed

Consequently, the foliations $\F_t$, $ \F _{\tau t}$ and $\F_{\tau^2 t}$ are projectively equivalent, for any $t$. 

For $t=1$, the projective equivalence of $\F_1, \F_\tau, \F_{\tau^2}$ was mentioned  in the begining of  Section \ref{ellipticpencilssurfaces}.  

Another example, the elliptic pencils \[\F_{6-3\tau} = [195, 49, 76, 67,3]\]  
and
\[\F_{\tau (6-3\tau)} =  \F_{3+9\tau} = [195, 67, 49, 76, 3]\]

\begin{pro}

The order four projective automorphism
\[T_{1,\infty}(x: y: z)) = (x+y+z : x + \tau^2 y + \tau z: x + \tau y + \tau^2 z)\]
acts on the parameter $t$ of $\F_ t$ as the involution 
\[t  \mapsto \frac{t+2}{t-1} =: \mbox{inv}_T(t),\]
in particular $T_{1,\infty}(\F_1) = \F_\infty$. It has the effects  
\[T_{1,\infty}(P_3(1)) = P_3(\infty),\quad T_{1,\infty} (P_3 (\infty)= P_3(1))\]
\[T_{1,\infty}(P_3(\tau )) = P_3(\tau^2),\quad T_{1,\infty}(P_2 (\tau^2)= P_3(\tau))\]
but is not an involution in the plane. 
\end{pro}

Proof:

The automorphism is represented by the matrix 
\[M = 
\begin{vmatrix}
 1  & 1 &1  \\
1   & \tau^2  & \tau \\
1 & \tau & \tau^2 \\
\end{vmatrix}
\]
for which  $M^4= \lambda \cdot \mbox{Id}$. It is not an involution in the plane, in fact,  
\[\begin{cases} T_{1,\infty}(1:1:1) = (1:0:0),\quad T^2_{1,\infty}(1:1:1)= (1:1:1)\\
T_{1,\infty}(1:\tau:\tau^2) = (0:1:0),\quad T^2_{1,\infty}(1:\tau:\tau^2)= (1:\tau^2:\tau)\\ 
T_{1,\infty}(1:\tau^2:\tau) = (0:0:1),\quad T^2_{1,\infty}(1:\tau^2:\tau)= (1:\tau:\tau^2) 
\end{cases}\]
For finding  the effect on  the parameter  $t \mapsto  \frac{t+2}{t-1}$, is just a matter of taking pullback 
\[T_{1,\infty}^*(\Omega + t E )\]
and comparing with $\Omega + t E$. 
\qed

Recalling the involutions on the parameter associated to $Q_1,Q_\tau, Q_{\tau^2},Q_\infty$ 
\[\begin{cases} q_1  (t) = - t -1 \\
q_\tau(t)= -t - \tau\\
q_{\tau^2}(t) = -t - \tau^2 = -t + 1 +\tau \\ q_{\infty}(t) = \frac{1}{t}\end{cases}\]
and the involution  $\mbox{inv}_T(t) = \frac{t+2}{t-1}$,  we remark that 
\[\mbox{inv}_T \circ q_1 \circ \mbox{inv}_T = q_\infty\]

Our method of realization of  elliptic pencils $\F_t$ extends directly  to  any parameter $t \in \mathbb{Q}(\tau)$  reached after  a finite  number  of compositions  of $q_1, q_\tau, q_{\tau^2}, q_\infty, \mbox{inv}_T$   applied to an element  of $\mathbb{Z}(\tau)$.

But we do not know if it is possible  to reach in algorithmically way all parameters  in  $\mathbb{Q}(\tau)$.

At last, a word on the involution $- t$ of the parameter.

Initially it was  a surprising experiment the diference between the elliptic pencils associated to $t$ and  to $-t$. For example, 
\[\begin{cases}\F_{-5 + 12 \tau}: \,  [231, 84, 67, 79, 1]\\
\F_{5 - 12 \tau}:\, [693, 208, 259, 223, 3]
\end{cases} \]
or 
\[\begin{cases}
\F_{1 + 13 \tau} : \, [477, 169, 133, 172, 3]\\
\F_{-1 - 13 \tau} : [159, 49, 61, 48, 1]
\end{cases}\]
But we remarked  that in the first example, for $t = - 5 +12 \tau$, while $-5 +12 \equiv 1  \,  \mbox{mod} (3)$,  it holds $5 -12 \equiv 2 \, \mbox{mod} (3)$. 
In the second,  while  $ 1 + 13 \equiv 2 \,  \mbox{mod}  (3)$,  it holds  $-1-13 \equiv  1  \,  \mbox{mod} (3)$.  

These arithmetic changes  and theirs  effects on degrees of the generic element of the pencils are compatible with the degree formula of  \cite{Puchuri}.

{Liliana Puchuri, Pontificia Universidad Cat\'olica del Per\'u, Per\'u.

email: lpuchuri@pucp.pe

Lu\'{i}s  Gustavo Mendes, Universidade Federal do Rio Grande do Sul, UFRGS, Brazil. 

email:   gustavo.mendes@ufrgs.br}

\end{document}